\documentclass{amsart}

\usepackage{amssymb}

\usepackage[all]{xy}

\newtheorem{thm}{Theorem}

\newtheorem{lem}[thm]{Lemma}
\newtheorem{cor}[thm]{Corollary}

\newtheorem{prop}[thm]{Proposition}

\theoremstyle{definition}
\newtheorem{defn}[thm]{Definition}

\newtheorem{say}[thm]{}

\newtheorem{const}[thm]{Construction}   

\newtheorem{rem}[thm]{Remark}          

\newtheorem*{ack}{Acknowledgments}      

\newtheorem{defn-thm}[thm]{Definition--Theorem}  
\newtheorem{defn-lem}[thm]{Definition--Lemma}  

\theoremstyle{remark}

\let \cedilla =\c
\renewcommand{\c}[0]{{\mathbb C}}  

\renewcommand{\o}[0]{{\mathcal O}} 
\newcommand{\z}[0]{{\mathbb Z}}

\renewcommand{\a}[0]{{\mathbb A}}

\newcommand{\p}[0]{{\mathbb P}}
\newcommand{\f}[0]{{\mathbb F}}

\newcommand{\map}[0]{\dasharrow}
\newcommand{\qtq}[1]{\quad\mbox{#1}\quad}

\newcommand{\proj}[0]{\operatorname{Proj}}

\newcommand{\aut}[0]{\operatorname{Aut}}

\newcommand{\chr}[0]{\operatorname{char}}

\newcommand{\tsum}[0]{\textstyle{\sum}}

\def\into{\DOTSB\lhook\joinrel\to}

\def\loccoh#1.#2.#3.#4.{H^{#1}_{#2}(#3,#4)}

\DeclareMathAlphabet{\mathchanc}{OT1}{pzc}%
                                {m}{it}

\newcommand{\sym}[0]{\operatorname{Sym}}
\newcommand{\hgt}[0]{\operatorname{ht}}

\begin{document}
\bibliographystyle{amsalpha}


\title[Families of elliptic curves]{Quadratic families of elliptic curves and\\
   unirationality of degree 1 conic bundles} 
\author{J\'anos Koll\'ar and Massimiliano Mella}

\maketitle

Let $K$ be a number field and $a_i(t)\in K[t]$ polynomials of degree 2.
We consider the family of elliptic curves
$$
E_t:=\bigl(y^2=a_3(t)x^3+a_2(t)x^2+a_1(t)x+ a_0(t)\bigr)\subset \a^2_{xy}
\eqno{(*)}
$$
parametrized by $t\in K$. Our aim is to show that there are many values
$t\in K$ for which the corresponding  elliptic curve $E_t$ has rank $\geq 1$.
Conjecturally this should hold for a positive proportion of them;
see the survey  \cite{MR1920278}. 
 We prove  rank $\geq 1$  for about the square root of all
$t\in K$, listed by height.

We are only interested in {\it nontrivial}  families, when
at least two of the curves $E_t$ are  smooth, elliptic and not
isomorphic to each other over $K$.
Thus   $a_3(t)$ is not identically 0 and not 
all the $a_i(t)$  are
constant  multiples of the same square $(t-c)^2$.

We  view the whole family as a single  algebraic surface in $\a^3_{xyt} $
and look at the distribution of $K$-points. The resulting surface
 has  degree 5 but its closure in $\p^3$ is very singular.
We prove the following.

\begin{thm}\label{main.ell.thm}  
Let $k$ be any field of characteristic $\neq 2$ and
$a_0(t),\dots, a_3(t)\in k[t]$ polynomials of degree 2
giving a nontrivial family of elliptic curves.
  Then the  surface
$$
S:=\bigl(y^2=a_3(t)x^3+a_2(t)x^2+a_1(t)x+ a_0(t)\bigr)\subset \a^3_{xyt}
\eqno{(**)}
$$
 is unirational over $k$.
\end{thm}

The proof has two parts. First assume that we  know a $k$-point
$p\in S$ that is not a 6-torsion point on the corresponding elliptic curve.
Then we have geometrically clear and  quite explicit formulas to 
 prove unirationality. The second, harder part is to show that there
are such $k$-points. We have not been able to turn this part of the proof
into explicit formulas; see Remark \ref{18.8.rem}.

From any sufficiently general $k$-point of $S$
we obtain families of   elliptic curves of rank $\geq 1$.

\begin{cor} \label{main.ell.cor} 
Let $k$ be an infinite  field of characteristic $\neq 2$. Then there
 are infinitely many different 
rational functions
$q(u)\in k(u)$ that are quotients of degree 2 polynomials
such that the rank of  the elliptic curve  $E_t$ as in $(*)$
is $\geq 1$ for all but finitely many values  $t=q(u)$ where $u\in k$. 
\end{cor}

Unirationality can also be used to exhibit points of small height.
For a Zariski open subset $U\subset S$, let 
 $N(U,B)$ be the number  of $K$-points of height $\leq B$ in $U$.
 Manin's  conjecture \cite{MR974910}  
 suggests that  $N(U,B)$
should grow  at least  like $B^{3/2}$, 
using the  naive height function  $\hgt(x,y,t):=\hgt(x)+\hgt(t)$.
 Theorem  \ref{main.ell.thm} 
implies that $N(U,B)$ grows at least like a 
power of $B$.

\begin{cor} \label{main.cor.2} Let $K$ be a number field. 
There is an $\epsilon>0$ such that for every $S$ as above
and for every Zariski open subset $U\subset S$,
$$
N(U,B)\geq c(S)\cdot B^{\epsilon} \qtq{for} B\gg 1,
$$
where $c(S)>0$  depends on $S$.
\end{cor}

The proof   gives an  explicit value for $\epsilon$
but it seems to be small.

\begin{say}[Connection with conic bundles]\label{spec.cases.say}
We can rewrite 
$$
a_3(t)x^3+\cdots + a_0(t)=A(x)t^2+B(x)t+C(x)
$$
where $A,B,C$ are cubics. Projection to the $x$-axis exhibits $S$ as 
the family of conics
$$
F_x:=\bigl(y^2=A(x)t^2+B(x)t+C(x)\bigr)\subset \a^2_{yt}.
$$
The conic $F_x$ is singular iff $x$ is a  root of the discriminant
$B(x)^2-4A(x)C(x)$; in general this happens 
for 6 different values of $x$. We get a possible 7th singular
fiber at infinity. After suitable birational transformations
the fiber at infinity is isomorphic to 
$$
F_{\infty}:=\bigl(\tilde a_3(s,t)=0\bigr)\subset \p^2_{stw},
$$
where $\tilde a_3(s,t):=s^2a_3(t/s)$ is the homogenization of $a_3(t)$. 
Thus  $S$ is birational to a conic bundle with $\leq  7$ singular fibers; 
see Definition \ref{CB.defns}. 

A very degenerate case is when $B(x)^2-4A(x)C(x)\equiv 0$. These are easy to 
enumerate by hand. The only non-rational surface occurs 
when $a_i(t)=c_i(t-\alpha)^2$ for every $i$. Then
setting $z:=y/( t-\alpha)$ gives the new equation
$$
z^2=c_3x^3+c_2x^2+c_1x+ c_0.
$$
Thus $S$ is birational to the product of an elliptic curve with $\p^1$.

 Our main technical result, Theorem \ref{main.thm},  says that
every  conic bundle   with 7 singular fibers  
is  unirational over its field of definition.
Not every conic bundle   with 7 singular fibers can be written as $(**)$
but unirationality is easier for the other ones; see
Section \ref{Uni.G.E.SEC}.


 Various coincidences among the $a_i(t)$
lead to simpler conic bundles. 
If $a_3(t)$ is reducible in $k[t]$ then $F_{\infty}$
is reducible over $k$ and we can contract either of its 
irreducible components to get a  conic bundle with 6 singular fibers.

If  $a_i(t)=c_iq(t)$ for every $i$ 
then the substitution  $y= zq(t)$ transforms the equation into
$$
q(t)z^2=c_3x^3+c_2x^2+c_1x+ c_0.
$$
The corresponding surface is  a  conic bundle with 4 singular fibers.
In this special case   Manin's  conjecture \cite{MR974910}  
 suggests that  $N(S,B)$
should grow   like $B^{2}$. 
In the complex multiplication case 
these were studied in \cite{MR2792990} where the bound
$N(S,B)\geq \epsilon\cdot B^{2-\eta}$ is proved for every $\eta>0$.
\end{say}

\begin{say}[Connection with higher dimensional conic bundles]
Let $k$ be an algebraically closed field. A (birational) conic
bundle over $k$ is a morphism  $\pi:X\to Y$ whose generic fiber is
isomorphic to a conic. It is a long standing open  problem to understand
which conic bundles are rational or unirational.
 These questions
are interesting only when $Y$ itself is unirational.

In many cases, for instance when $\dim X=3$, one can realize
$Y$ as a conic bundle  $\tau:Y\map B$, hence we have
$\tau\circ\pi:X\map B$ whose generic fiber $X_B$ is a 2-dimensional
conic bundle over the function field $k(B)$.
Thus if   $X_B$ is (uni)rational over $k(B)$ then $X$ is (uni)rational over $k$.

The conjecture of  \cite{MR899398} says that rationality is
almost equivalent to $\bigl(K_{S_B}^2\bigr)\geq 5$.
Our Corollary \ref{main.thm.2} implies that  
 $X$ is unirational if $\bigl(K_{S_B}^2\bigr)\geq 1$.
\end{say}

\begin{ack} We
 thank J.-L.~Colliot-Th\'el\`ene, Y.~Liu, R.~Munshi, C.~Skinner, 
Y.~Tschinkel for
comments, discussions and references and  S.~Kov\'acs for his
hospitality at the University of Washington where this project started.
We received numerous corrections and suggestions from the referee.

Partial financial support  to JK  was provided  by  the NSF under grant number
 DMS-1362960 and to MM by  {\it Geometria sulle Variet\`a Algebriche} (MUIR).
\end{ack}

\section{Minimal conic bundles}

\begin{defn}[Conic bundles]\label{CB.defns}
 Let $k$ be a field. 
A surface with a {\it pencil of rational curves} 
over $k$ is a  projective, geometrically irreducible
surface $T$ together with a  morphism  $\pi:T\to B$ 
to a smooth, projective curve $B$
such that  the geometric generic fiber of $\pi$
is isomorphic to $\p^1$. M.~Noether proved in 1870 that if  $k=\c$
then  $T$ is birational (over $B$) to the 
trivial family  $\p^1\times B\to B$  \cite{MR1509694}.

This is no longer true if $k$ is not algebraically closed
and the birational properties of such surfaces can be quite subtle.
The generic fiber  $F_{k(B)}$ is isomorphic to a conic
and $T$ is birational (over $B$) to  $\p^1\times B$
iff  $F_{k(B)}\cong \p^1_{k(B)}$. The latter holds iff
 $F_{k(B)}$ has a $k(B)$-point. This in turn is equivalent to  $\pi:T\to B$
having a section.

A {\it  conic bundle}
over $k$ is a smooth, projective, geometrically irreducible
surface $S$ together with a  morphism  $\pi:S\to B$ 
to a smooth, projective curve $B$
such that  geometric fibers of $\pi$ are plane 
conics (either smooth or a  pair of lines).

 A  conic bundle is called {\it  minimal} if it can not be
obtained from another  conic bundle  $\pi_1:S_1\to B$
by blowing up points. 
If every fiber of $\pi$ is smooth
then $\pi:S\to B$ is minimal; these are the trivial examples.

As a generalization of Noether's result, Iskovskikh \cite{MR525940}
and  Mori \cite{Mori82} proved that every surface with a 
 pencil of rational curves  is birational (over $B$) to a 
minimal conic bundle. Thus, up-to birational equivalence,
it is sufficient to study minimal conic bundles. 
\end{defn}

It is a long standing interesting question to understand which
conic bundles are rational or unirational.
If   $\pi:S\to B$ is unirational then there is a dominant
morphism $\p^2\map B$, hence  $B\cong \p^1$ by L\"uroth's theorem;
see also Lemma \ref{unirtl.iff.multisecton}.
Thus from now on we restrict our attention to
minimal conic bundles $\pi:S\to \p^1$.

A basic numerical invariant of a minimal conic bundle is the
{\it number of singular fibers,} denoted by $\delta(S)$.
People coming from the theory of del~Pezzo surfaces prefer to use instead
the self intersection of the canonical class  $(K_S^2)$,
also called the {\it degree}. These are related by the formula
$$
(K_S^2)=8-\delta(S).
$$
As with  del~Pezzo surfaces,
the birational complexity of $S$ increases 
as the degree decreases. 

Conic bundles with  $\delta(S)\in\{0,1,2,3\} $  are rational
if they have a $k$-point; see \cite[Sec.IV.8]{MR0460349}
 or \cite[III.3.13]{rc-book}.

Conic bundles with  $\delta(S)\in\{4,5\} $
 are unirational if they have a $k$-point. 
This has been essentially known to Segre \cite{segre51},
Iskovskikh  \cite{MR0271109} and  Manin  \cite{MR0460349},
 though most of their work
focused on smooth del~Pezzo surfaces over infinite fields.
Extra  difficulties over small fields were resolved in \cite{MR1956057}. 
For the conic bundle case see also \cite{MR1383424}.

Unirationality of  degree 2  del~Pezzo surfaces with
a $k$-point  is still not fully proved;
see \cite{MR3245139, 2014arXiv1408.0269F} for the known results.

We treat   
conic bundles with $\delta(S)=6$ in Section \ref{deg2.section}.
\medskip

In this note we settle the next case and prove the following

\begin{thm}\label{main.thm} Let $k$ be a field of characteristic $\neq 2$
and $\pi:S\to \p^1$ a 
  conic bundle with 7 singular fibers. Then $S$ is unirational over $k$.

Moreover, if $|k|\geq  53$ then $\pi$ has a degree 8 rational multi-section. That is,
there is a map  $\sigma:\p^1\to S$ such that
$\pi\circ \sigma:\p^1\to \p^1$ has degree 8.
\end{thm}


Combining with the earlier results gives the following 
more uniform statement.

\begin{cor}\label{main.thm.2} Let $k$ be a field of characteristic $\neq 2$
and $\pi:S\to \p^1$ a conic bundle with $\leq 7$ singular fibers. 
The following are equivalent.  
\begin{enumerate}
\item  $S$ has a $k$-point.
\item $S$  is unirational over $k$.
\end{enumerate}
\end{cor}

Conic bundles with $\geq 9$ singular fibers seem to behave quite
differently and various heuristics suggest contradictory possibilities.
Working with forms similar to $(**)$ shows that such conic bundles
appear as quadratic families of genus $\geq 2$ hyperelliptic curves.
This suggests that they should contain  few $k$-points
and they should not be unirational.  On the other hand, they 
become rational after a finite field extension, so they are not far from
being unirational. 

This leaves conic bundles with $8$ singular fibers as a quite interesting
case. Some of these are unirational by \cite{MR1309090}.
 We hope to return to them in the future.

Unirationality of conic bundles---with any number of singular fibers---over 
local fields and sufficiently large finite fields
  is treated in
\cite{MR826395, MR1381777, k-loc, MR2019976}.
 These, however, seem to be   special properties
of these fields rather than  of conic bundles.

\begin{say}[Proof of Corollary \ref{main.ell.cor}]
Pick any rational point $p\in S(k)$ that lies on a  smooth fiber $F$.
Then $F\cong \p^1_k$ and projection of $F$ to the $t$-axis is described
by a degree 2 rational function $q(u)$. We are done unless 
we get torsion points for infinitely many  $u\in K$. 
By \cite[Thm.III.11.4]{sil-adv}  
this can happen only if the whole fiber $F$ consists of $m$-torsion points
for some fixed $m$, in which case
 the modular curve  $X_1(m)$ is rational. Thus 
$m\leq 12$ and there are only  finitely many such fibers. \qed
\end{say}

\begin{say}[Proof of Corollary \ref{main.cor.2}]
A map $\p^2\map S\map \p^1_x\times \p^1_t$ is described by
4 homogeneous polynomials of some degree $d$:
$$
(u{:}v{:}w)\mapsto  \bigl(H_1(u{:}v{:}w):H_2(u{:}v{:}w), 
H_3(u{:}v{:}w):H_4(u{:}v{:}w)\bigr).
$$
Thus if $\hgt(u{:}v{:}w)\leq B^{1/d}$ then the image point has
height $\leq B$  (up to constant factors). This gives $\epsilon\geq \frac2{d}$.
\qed
\end{say}

A difficulty in computing the best $\epsilon$ is that the map
$\phi:\p^2\map S$ does not determine $d$ since composing $\phi$
with a birational self-map of $\p^2$ can increase or decrease 
the value of $d$.
The geometric description suggests that it might be possible to find
$H_i$ of degree $\leq 8$, but so far we have managed to write down
only much higher degree examples.

\section{The geometry of conic bundles}

We give a quick summary of the geometry of rational surfaces and conic bundles.
This is a classical topic. For some modern treatments see
\cite[Secs.III.2--3]{rc-book},   \cite[Chap.3]{ksc} or \cite[Chap.8]{dolg-cl}.

\begin{say}[Rational surfaces] 
Let $\bar k$ be an algebraically closed field and $S$ a smooth 
rational surface  over $\bar k$. Set  $r:=(K_S^2)$.

Basic examples are $\p^2$ and 
$\f_e:=\proj_{\p^1}\bigl(\o_{\p^1}\oplus \o_{\p^1}(e)\bigr)$
for  $e\geq 0$. 

Every other $S$ is obtained either from $\p^2$ by blowing up $9-r$ points
or from some $\f_e$   by blowing up $8-r$ points.
(There can be many such representations.)

Note that  $h^0(\p^2, -K)=10$,
$h^0(\f_e, -K)=9$ for $0\leq e\leq 3$ and 
$h^0(\f_e, -K)=e+6$ for $e\geq 3$. Thus we see that
$\dim |-K_S|\geq r$.

If $C\in |-K_S|$ is irreducible and reduced then, by the
adjunction formula,   
$$
2p_a(C)-2=\bigl(C\cdot (C+K_S)\bigr)=0, \qtq{hence} p_a(C)=1.
$$
Otherwise every irreducible component of $C$ is a smooth rational
curve.
\end{say}

\begin{defn}[Weak del~Pezzo surfaces]\label{wdP.all.say}
Classically, a smooth, projective, geometrically irreducible surface 
$S$ is called a {\it  del~Pezzo} surface if  $-K_S$ is ample, but modern usage
 allows $S$ to have Du~Val singularities.

A smooth, projective, geometrically irreducible surface
$S$ is called a {\it weak del~Pezzo} surface if $(K_S^2)>0 $
and the following
equivalent conditions hold.
\begin{enumerate}
\item $S$ is the minimal resolution of a del~Pezzo surface  with
Du~Val singularities.
\item $-K_S$ is nef, that is, $(C\cdot K_S)\leq 0$ for every 
irreducible curve $C\subset S$.
\item $|-K_S|$ has no fixed components.
\item There is an irreducible curve  $C\in |-K_S|$.
\end{enumerate}
In these cases $\dim |-K_S|= (K_S^2) $
is called the {\it degree} of $S$. 
(If   $-K_S$ is ample and 
$d:=(K_S^2)\geq 3$ then  $|-K_S| $ gives an embedding $\phi_S:S\into \p^d$
as a surface of degree $d$.)  For a curve  $C\subset S$ the
intersection number  $(C\cdot K_S)$ is called the
 {\it degree} of $C$. (If $d:=(K_S^2)\geq 3$ then this is the same as the
degree of $\phi_S(C)\subset \p^d$.)
\end{defn}

\begin{defn}[Minimal conic bundles]  For any conic bundle  $\pi:S\to \p^1$ we 
use $F$ to denote a fiber of $\pi$. The adjunction formula gives that
$(F\cdot K_S)=-2$.

Assume next that  $\pi:S\to \p^1$ is minimal. Let $p\in \p^1$ be any point and 
$F_p$ the fiber of $\pi$ over $p$. 
Then a fiber $F$ is isomorphic, over $k(p)$,  either to a smooth conic or 
to a  conjugate pair of lines. In particular, $F_p$ is irreducible
 over $k(p)$. 

If  $\pi:S\to \p^1$ is a degree 1, weak del~Pezzo conic bundle
then the linear system  $|-K_S|$
is a pencil of elliptic curves with a unique base point; henceforth denoted
by $p^*$. We use $F^*$ to denote the fiber containing $p^*$.
\end{defn}

\begin{defn}[Bertini involution]\label{bert.inv.defn} 
(See \cite[Sec.8.8]{dolg-cl} or 
\cite[3.42]{ksc}.)
Let $S$  be a  del~Pezzo surface of degree 1, 
possibly with Du~Val singularities. Then 
  $|-2K_S|$ gives  a morphism  of degree two
$\pi_0:S\to Q$  (where $Q$ is a quadric cone in $\p^3$).
The base point $p^*$ is the preimage of the vertex of $Q$.
 The corresponding Galois involution $\tau:S\to S$ is called the 
{\it Bertini involution.}

Let $C\in |-K_S|$ be a geometrically  irreducible curve;
it has arithmetic genus 1.
Using the base point $p^*\in C$ as the identity of the group structure,
the Bertini involution sends  $q\in C$ to $-q\in C$.
\end{defn}

\begin{say}[Curves of low degree] \label{low.deg.curves}
Let $T$ be a degree 1 del~Pezzo surface, possibly with Du~Val singularities.
The following is the list of degree $1$ curves on $T$;
with degree as in Definition \ref{wdP.all.say}.
\begin{enumerate}
\item Double covers of a ruling of $Q$. These are members of $|-K_T|$.
\item In some cases $\pi_0^{-1}(Q\cap P)$ is reducible for a plane
$P\subset \p^3$ not passing through the vertex of $Q$.
 (If $T$ is  smooth,  this happens for 120 planes.)
Both components
are then smooth, rational. These are the $(-1)$-curves on $T$; 
they do not pass through $p^*$.
\end{enumerate}
The following is the list of degree $2$ curves on $T$.
\begin{enumerate}\setcounter{enumi}{2}
\item Usually  $\pi_0^{-1}(Q\cap P)$ is irreducible for a plane
$P\subset \p^3$. Such curves 
do not pass through $p^*$.
\item In some cases $\pi_0^{-1}(Q\cap Q')$ is reducible for a quadric
$Q'\subset \p^3$.  
\end{enumerate}
\end{say}

 From the above list we can read off the connection between the 
Bertini involution and the conic bundle structure.

\begin{say}[Bertini involution and the conic bundle structure]
 \label{bertini.and.CB}
Let $T$ be a degree 1 del~Pezzo surface, possibly   with Du~Val singularities,
and $|F|$  a pencil 
of degree 2 rational curves.  
It has at least 1 member   $F^*$ passing through $p^*$.

If $F^*$ is reducible then at least one of the 
  irreducible components must be as 
in case (\ref{low.deg.curves}.1). 
Furthermore, either both  irreducible components are as 
in case (\ref{low.deg.curves}.1) or one of them is as 
in case (\ref{low.deg.curves}.1) and the other as in  case (\ref{low.deg.curves}.2). In the first case $p^*$ is the singular point of $F^*$ and 
 $F^*$ is  invariant under the Bertini involution.
In the second case $p^*$ is a smooth point of $F^*$ and 
 $F^*$ is not  invariant under the Bertini involution.

If $F^*$ is  irreducible
 then $F^*$  must be as in case (\ref{low.deg.curves}.4).
Thus $F^*$ is not invariant under the Bertini involution.

Note further that $\pi_0^*|\o_Q(1)|$ is a complete linear system,
thus if one member of  $|F|$  is  invariant under the Bertini involution
then every  member is  invariant. 

Furthermore, if a pencil  $|F|$  is  invariant under an involution
then at least 2 members are invariant, hence, in our case,
every  member is  invariant. 

\end{say}

We have proved the following.

\begin{prop}  \label{bertini.and.CB.prop}
Let $\pi:S\to \p^1$ be a  weak del~Pezzo conic bundle of degree 1.
Then the Bertini involution preserves the  conic bundle structure iff
$p^*$ is a singular point of  a singular fiber of $\pi$.\qed
\end{prop}

\section{Classification up-to isomorphism}

\begin{lem}\label{irred.on.CB.lem}
Let $\pi:S\to B$ be a  minimal conic bundle with at least 1
singular fiber.  Then 
 every line bundle on $S$ is of the form
$\o_S(-aK_S)\otimes \pi^*L_B$ where $L_B$ is a line bundle on $B$.
In particular,  every curve $C\subset S$ is either a fiber or
 the projection $\pi|_C:C\to B$
has even degree.  
\end{lem}

Proof. Assume for simplicity that $k$ is perfect.
Let  $L$ be a line bundle on $S$ and $F\subset S$ a reducible geometric fiber
with  irreducible components $F'+F''$. 
If $(L\cdot F')\neq (L\cdot F'')$ then the conjugates of $F'$  form
a conjugation invariant set of pairwise disjoint $(-1)$-curves.
This contradicts the minimality assumption.  Thus
$(L\cdot F')= (L\cdot F'')$ hence 
 $L$ has even degree on the generic fiber. Therefore
$L(aK_S)$ is trivial on every smooth fiber of $\pi$ for some $a\in\z$.
By the above considerations, 
 $L(aK_S)$ is also trivial on the singular fibers. By 
cohomology and base change, this implies that 
$$
L(aK_S)\cong \pi^*\pi_*\bigl(L(aK_S)\bigr).
$$
Essentially the same argument works if $k$ is not perfect; see
\cite[Chap.2]{Mori82}. \qed

\begin{lem}\label{exceptional.cases.prop}
Let $k$ be a field, $B$ a smooth conic and 
  $\pi:S\to B$  a minimal  conic bundle of degree $r:=(K_S^2)\geq 1$. 
Then one of the following holds.
\begin{enumerate}
\item $S$ is a  weak del~Pezzo surface,
\item  $r=2$ and $|-K_S|=C+|2F|$ where $C$ is a conjugate pair of
disjoint smooth rational curves with self intersection $-3$.
\item  $r=1$ and $|-K_S|=C+|F|$ where $C$ is a geometrically 
irreducible  smooth rational curve with self intersection $-3$.
\item  $\pi:S\to \p^1$ is a $\p^1$-bundle.
\end{enumerate}
\end{lem}

Proof. If $\pi$ has a section then $\pi:S\to \p^1$ is a $\p^1$-bundle.
Thus assume from now on that there are no sections.

Write   $|-K_S|=C+|M|$ where $C$ is the fixed part and
$|M|$ the mobile part.  If $C=0$ then $S$ is a  weak del~Pezzo surface;
see Definition \ref{wdP.all.say}.3.

If $C\neq 0$ then 
$C$ can not contain a fiber, so $C$ is a union of  
multi-sections. Furthermore,
 $$
2=(-K_S\cdot F)=(C\cdot F)+(M\cdot F)\geq (C\cdot F)\geq 2,
$$
the last inequality by Lemma \ref{irred.on.CB.lem}.
Thus $(C\cdot F)=2$ and so we have 3 possibilities for 
$C$:  geometrically irreducible,   conjugate  pair of sections or geometrically nonreduced; the latter can happen only  in characteristic 2.  
In all cases $(M\cdot F) =0$, thus $M\sim sF$
 for some $s$ and 
 $s=\dim |sF|= \dim |M|=\dim |-K_S|\geq r$.
Note that  $\deg \omega_C=\bigl(C\cdot (C+K_S)\bigr)=-s(C\cdot F)=-2s$
and $(C^2)=\bigl((-K_S-sF)^2\bigr)=r-4s$.
We distinguish several cases.

If $C$ is  geometrically irreducible and reduced then $\deg \omega_C\geq -2$. 
Thus $s=1$ hence $r=1$ and $C$ is a smooth, geometrically rational curve
with  $(C^2)=-3$; this is case (3).

 If $C$ is  disconnected then  $C_{\bar k}$ is the disjoint union of 2 components that are conjugate over $k$. Thus 
either  $\deg \omega_C\geq 0$ (this is not possible for us) or 
$\deg \omega_C= -4$ and $C_{\bar k}$ is the disjoint union of 2 
smooth,  rational curves. If $r=s=2$ then  $(C^2)=-6$; this is case (2).
Otherwise $r=1$ and $(C^2)$ is odd, but this is not possible.

If $C$ is  connected, geometrically reduced and reducible
then $\deg \omega_C\geq -2$, hence $s=r=1$ and $(C^2)=-3$. 
In this case  $C_{\bar k}=C_1+C_2$ and so $(C^2)=2(C_1^2)+2(C_1\cdot C_2)$ is even, 
a contradiction. This case can not happen.

Finally, if  $C$ is  irreducible but geometrically nonreduced
then $C$ intersects only 1 of the geometric irreducible components of every
singular fiber. This is impossible by Lemma \ref{irred.on.CB.lem}.
 \qed

\begin{cor} \label{pts.exist.cor}
A degree 1 conic bundle over $k$ has a $k$-point.
\end{cor}

Proof. If $S$ is weak del~Pezzo then the unique base point of
$|-K_S|$ is a $k$-point. If $S$ is as in  (\ref{exceptional.cases.prop}.3)
then $C$ is a smooth, rational curve and $(C^2)=-3$ gives a degree 3 point on
 $C$. Thus $C\cong \p^1$. 

We still need to deal with the cases when $S$ is not minimal, thus
$S$ is obtained from a minimal conic bundle $T$ of degree $m+1$
by blowing up $m$ points. In particular, $T$ has a 0-cycle of degree 1.
For $m\geq 3$ this easily implies that $T$ has a $k$-point 
\cite[Sec.IV.8]{MR0460349}.

If $m=2$ then   $T$ is birational to a cubic surface with a 
degree 2 point, hence it has a $k$-point as well.

If $m=1$ then $T$ has a $k$-point (the one we blow up)
and so does $S$. \qed

\begin{rem} It is not known if a cubic surface with a 0-cycle of degree 1
has a $k$-point or not, see \cite{MR0429731}. There are degree 2
del~Pezzo surfaces with a degree 3 point---and hence with a 
0-cycle of degree 1---but without $k$-points.
\end{rem}

\begin{defn} We divide  degree 1, smooth,  minimal conic bundles 
$\pi:S\to \p^1$ into three types
up-to isomorphism. By Lemma \ref{exceptional.cases.prop}, either 
we are in case (\ref{exceptional.cases.prop}.3) or $|-K_S|$ is a pencil with a unique base point $p^*\in S(k)$. In the latter case
let  $F^*$ denote the fiber  of $\pi$ containing $p^*$.
We have three possibilities for $S$. 
\begin{enumerate}
\item[Type G:] (General case) $S$ is a weak del~Pezzo surface and $p^*$
 lies   on a smooth fiber  $F^*$. Therefore  $F^*\cong \p^1_k$.
\item[Type S:] (Special case) $S$ is a weak del~Pezzo  surface and $p^*$
 lies  on a singular fiber  $F^*$. In this case $p^*$ 
 is the unique singular point of $F^*$ .
\item[Type E:] (Exceptional case) $S$ is not a weak del~Pezzo  surface.
\end{enumerate}
(The birational classification is different. As a consequence of
 Theorem \ref{main.thm} and Proposition \ref{typeS.bir.typeG}
we see that every minimal conic bundle of degree 1 is birational
to a minimal conic bundle of Type G.)
\end{defn}

Next we give more detailed descriptions of the two main Types.

\begin{say}[Type G]\label{type.G.say}
After blowing up $p^*$,  the  pencil $|-K_S|$  becomes   base point free
and gives  a morphism  $\pi_2:B^*S\to \p^1_{uv}$ whose general fiber is
a genus 1 curve (smooth if $\chr k\neq 2,3$).
We thus have a morphism 
$$
(\pi, \pi_2): B^*S\to \p^1_{xy}\times \p^1_{uv}.
$$
Since $(-K_S\cdot F)=2$, the morphism $(\pi, \pi_2) $
has degree 2. 
Assume that there is an irreducible  exceptional curve $D$.
Then $D$ is contained in a fiber of $\pi$, but every fiber of
$\pi$ is irreducible by  Lemma \ref{irred.on.CB.lem}. Thus
$(\pi, \pi_2) $ is finite.  If  $\chr k\neq 2$, let
  $C_S\subset\p^1_{xy}\times \p^1_{uv}$ be its
 branch curve. 
A general fiber $F$ of $\pi$ has genus 0 and 
$\pi_2$ restricts to a double cover $F\to \p^1_{uv}$. Thus the latter has 2 branch points. A general fiber $E$ of $\pi_2$ has genus 1 and 
$\pi$ restricts to a double cover $E\to \p^1_{xy}$. Thus the latter has 4 branch points. Therefore $C_S$ has  bidegree $(2,4)$. 
We can thus write the equation of $C_S$ as
$$
a_4(u,v)x^4+a_3(u,v)x^3y+a_2(u,v)x^2y^2+a_1(u,v)xy^3+ a_0(u,v)y^4=0,
\eqno{(\ref{type.G.say}.1)}
$$
where the $a_i$ are homogeneous of degree 2.

If  $\chr k\neq 2$ then a double cover of a smooth surface 
is smooth iff the branch curve is smooth, hence $C_S$ is smooth.

We will repeatedly use 2 observations connecting the $a_i$ to the
properties of $B^*S$. 
\begin{enumerate}\setcounter{enumi}{1}
\item The curve defined by (\ref{type.G.say}.1) is singular along
$y=0$ iff there is a $(u_0, v_0)$ that is a double root of $a_4$ and root of $a_3$. In particular, if $a_4\equiv 0$ then $C_S$ is singular at the roots of $a_3$. 
\item The preimage of $(y=0)$ in $B^*S$ is reducible (resp.\ geometrically reducible) iff 
$a_4$ is a square in $k[u,v]$ (resp.\ in $\bar k[u,v]$).
\end{enumerate} 
We can choose the $x,y$-coordinates such that $F^*$ lies over $(y=0)$.
Note that the corresponding fiber of $\pi_B^*S\to \p^1_{xy}$ is the union of the 
exceptional curve $E^*\subset B^*S$ and of the
birational transform of  $F^*$.  Thus $a_4(u,v)$ is a square.
We can choose the $u,v$-coordinates such that $a_4(u,v)=v^2$.
Since $(y=v=0)$ is a smooth point of $C_S$, we see that $u^2$ appears in
$a_3(u,v) $ with non-zero coefficient.
With these choices, 
 in affine coordinates  $x,u,z$,   the equation of $B^*S$ is 
$$
z^2=x^4+a_3(u)x^3+a_2(u)x^2+a_1(u)x+ a_0(u)
\eqno{(\ref{type.G.say}.4)}
$$
 $\deg a_i\leq 2$ and $\deg a_3= 2$ where we use the shorthand $a_i(u):=a_i(u,1)$.

We can thus think of a conic bundle of Type G as a quadratic family
of genus 1 curves with two $k$-points.

 {\it Warning.} This double cover gives a 
biregular  involution
on $B^*S$, hence a
birational involution
on $S$ but, by   Proposition \ref{bertini.and.CB.prop},
 it is {\em not} the Bertini involution; see also
Paragraph \ref{explicit.say}.
\end{say}

\begin{say}[Type S]\label{type.S.say}
Here $p^*$ is the singular point of $F^*$ and, 
after blowing up $p^*$,  the  pencil $|-K_S|$ becomes   base point free.
We again have a morphism 
 $$
(\pi, \pi_2): B^*S\to \p^1_{xy}\times \p^1_{uv}.
$$
The birational transform of  $F^*$ 
(consisting of a conjugate pair of $(-2)$-curves) is contracted. 
As above,  if $\chr k\neq 2$ then, $(\pi, \pi_2) $
 branches along a  curve $B_S$ of bidegree $(2,4)$, however
$B_S$ has 2 singular points at the 2 images of the 
 birational transform of  $F^*$. These 2 points lie on the same fiber, namely $F^*$, of
the 1st projection. Thus $B_S=E_0\cup C_S$ where 
$C_S\subset\p^1_{xy}\times \p^1_{uv}$ 
is a smooth curve  of bidegree $(2,3)$.
We can thus write the equation of $C_S$ as
$$
a_3(u,v)x^3+a_2(u,v)x^2y+a_1(u,v)xy^2+ a_0(u,v)y^3=0,
\eqno{(\ref{type.S.say}.1)}
$$
where the $a_i$ are homogeneous of degree 2.
We can choose the $x,y$-coordinates such that $E_0=(y=0)$.
Then we have the following.
\begin{enumerate}\setcounter{enumi}{1}
\item  $a_3$ is irreducible in $k[u,v]$.
\item The curve defined by (\ref{type.S.say}.1) is singular along
$x=0$ iff there is a $(u_0, v_0)$ that is a double root of $a_0$ and root of $a_1$. 
\end{enumerate} 
Indeed,  the  conjugate pair of irreducible components of $F^*$ map to
the points  $(y=a_3=0)$, proving (2) while (3) is the same argument as
 (\ref{type.G.say}.2).

 In suitable coordinates 
we can write the affine equation as
$$
z^2=a_3(u)x^3+a_2(u)x^2+a_1(u)x+ a_0(u),
\eqno{(\ref{type.S.say}.4)}
$$
where $\deg a_i\leq 2$ and $\deg a_3= 2$. 

As before, this double cover gives a birational involution
on $S$. By Proposition \ref{bertini.and.CB.prop} it is  the Bertini involution.
\end{say}

\section{Unirationality for Type G and E}\label{Uni.G.E.SEC}

\begin{lem}[Enriques criterion] \label{unirtl.iff.multisecton}
A surface with a pencil of rational curves $\pi:T\to B$ is unirational iff
it has a rational multi-section. That is, a morphism
$\sigma:\p^1\to T$ such that $\pi\circ \sigma:\p^1\to B$ is nonconstant.
\end{lem}

Proof.  Assume that there is a dominant map  $p:\p^n\map T$. 
Over an infinite base field, a general line $L\subset \p^n$
gives a  morphism
$\sigma:\p^1\to T$ such that $\pi\circ \sigma:\p^1\to B$ is nonconstant.

Over finite fields one needs to be a little more careful; see
\cite[Lem.12]{MR1956057}. 

Conversely, assume that there is  a morphism
$\sigma:\p^1\to T$ such that $\pi\circ \sigma:\p^1\to B$ is nonconstant.
Then $\p^1\times_BT\to \p^1$ is a surface with a pencil of rational curves
that has a section. Thus, as we noted in
Definition \ref{CB.defns}, $\p^1\times_BT$ is rational 
hence $T$ is unirational.\qed

\begin{prop} A degree 1 conic bundle $\pi:S\to \p^1$ of Type E
 is unirational.
\end{prop}

Proof.  We are in case (\ref{exceptional.cases.prop}.3), thus
 $|-K_S|=C+|F|$ where $C$ is a geometrically 
irreducible  smooth rational curve with self intersection $-3$.
Thus $\o_S(C)|_C$ is a line bundle of odd degree, hence
$C\cong \p^1$. Thus $S$ is  unirational by (\ref{unirtl.iff.multisecton}). \qed

\begin{prop} \label{type.G.unirtl}
A degree 1 conic bundle $\pi:S\to \p^1$ of Type G
 is unirational.

More precisely,  the Bertini involution  (as in Definition \ref{bert.inv.defn})
does not preserve the conic bundle structure
and $\tau(F^*)\subset S$ is
 a  rational multi-section of degree $ 8$. 
\end{prop}

Proof. As in Paragraph \ref{bert.inv.defn}, $|-2K_S|$ gives a 
morphism  $\pi_0:S\to Q$ 
and, as we discussed in Paragraph \ref{bertini.and.CB}, 
 the conic bundle structure is given by
a pencil of rational curves $|F|$ coming from quadric sections  $Q'\cap Q$
whose preimages in $S$ are reducible. These preimages are thus members
of $|-4K_S|$ and the Bertini involution interchanges
$|F|$ with $|-4K_S-F|$.   Since  $\bigl(F\cdot (-4K_S-F)\bigr)=8$,
the curves in $|-4K_S-F|$ give degree 8 
multi-sections. 

The fiber $F^*$ is rational over $k$, thus $\tau(F^*)$ is
 a rational multi-section. So $S$ is unirational by the 
Enriques criterion \ref{unirtl.iff.multisecton}.
\qed

\begin{say}[Explicit formulas] \label{explicit.say}
Fixing a value of $u=u_0$ set $a_i:=a_i(u_0)$. We get an 
 elliptic curve
$$
z^2=x^4+a_3x^3+a_2x^2+a_1x+a_0.
$$
The base point $p^*$ is the point at infinity on the $(z\sim x^2)$-branch.
Then the linear system  $|2p^*|$ is given by
parabolas of the form
$$
z=x^2+\tfrac12 a_3x+\lambda.
$$
Note that
$$
\bigl(x^2+\tfrac12 a_3x+\lambda\bigr)^2=x^4+a_3x^3+a_2x^2+a_1x+a_0
$$
is a degree $\leq 2$ equation in $x$ and there is a unique choice, namely 
$$
\lambda= \tfrac18 (4a_2-a_3^2),
$$ 
that gives a linear equation. 
The root of this equation gives the $x$-coordinate of
$-p'$ where  $p'$ is the point
 at infinity on the $(z\sim -x^2)$-branch. Thus
$-p'$ is given by
$$
\Bigl(x:=\frac{(4a_2-a_3^2)^2-64 a_0}{64a_1-8a_3(4a_2-a_3^2)},
z:=x^2+\tfrac12 a_3x+\tfrac18(4a_2-a_3^2)\Bigr).
$$
Thus, in the conic bundle case, the restriction of
$\tau$ to $F^*$ is given by
$$
\begin{array}{rcl}
x(u)&=&\displaystyle{
\frac{(4a_2(u)-a_3(u)^2)^2-64 a_0(u)}{64a_1(u)-8a_3(u)(4a_2(u)-a_3(u)^2)}}
\qtq{and}\\[2ex]
z(u)&=&x(u)^2+\tfrac12 a_3(u)x(u)+\tfrac18(4a_2(u)-a_3(u)^2).
\end{array}
$$
In the type G case $\deg a_3=2$, so
 the numerator of $x(u)$ has degree 8 and the denominator
degree 6. In particular, $x(u)$  is never constant, in agreement with 
our earlier computations in Paragraph \ref{bertini.and.CB}.
\end{say}

\section{Unirationality for Type S; first construction}

In the special case there are no obvious rational curves on $S$
but the smooth fibers of $\pi$ are conics and usually
some of them are rational. However, the Bertini involution
maps every fiber to itself, so it does not yield new rational curves.

We look instead at a degree 4 endomorphism  $m_2:S\map S$
obtained as follows.
(More generally, we have a degree $r^2$ endomorphism  $m_r:S\map S$
for every $r\in \z$.)

\begin{say}[Doubling map] \label{doubling.map.say}
Let $\pi:S\to \p^1$ be a degree 1, minimal, weak del~Pezzo conic bundle
of type S.
Let $C\in |-K_S|$ be a geometrically  irreducible curve.
Using the base point $p^*\in C$ as the identity of the group structure,
the doubling map  $m_2$ sends  $q\in C$ to $2q\in C$.
Note that $q, -q$ lie on the same fiber and so
 $2q, -2q$ also lie on the same fiber. Thus the $x$-coordinate 
of $m_2(q)$ depends only on $\bigl(x(q),u(q)\bigr)$.

So pick a smooth fiber $F$, say over $x=0$. 
Let $q_0\in F(\bar k)$ have coordinates  $(0,y_0,u_0)$ .
By  explicit computation we get that
$$
x\bigl(m_2(q_0)\bigr)=\frac{a_1(u_0)^2-4a_2(u_0)a_0(u_0)}{4a_3(u_0)a_0(u_0)}.
\eqno{(\ref{doubling.map.say}.1)}
$$
We know that $a_3, a_0$ are not identically zero 
since the vanishing of either one would give a reducible, hence singular, branch curve $C_S$.
If this expression is non-constant then $m_2(F)\subset S$
is a multi-section of $\pi$. 
Next we describe the cases when this construction  does not give a 
rational multi-section.
\end{say}

\begin{prop} \label{2.4.specv.sec}
Let $\pi:S\to \p^1$ be a degree 1, minimal conic bundle 
of Type S and
 $F_0\subset S$  a smooth fiber of $\pi$.
Then   one of the following holds.
\begin{enumerate}
\item $m_2(F_0)\subset S$ is a multi-section of $\pi$,
\item $m_4(F_0)\subset S$ is a multi-section of $\pi$,
\item    $F_0$ consists of 2-torsion points and 
$m_r(F_0)\subset F_0\cup \{p^*\}$  for every $r\in\z$.
\item    $F_0$ consists of 3-torsion points and 
$m_r(F_0)\subset F_0\cup \{p^*\}$  for every $r\in\z$.
\end{enumerate}
\end{prop}

Proof. We are in case (1) unless (\ref{doubling.map.say}.1)
 is constant.  
Then, for some $\alpha\in k$,
$$
a_1^2-4a_2a_0=4\alpha a_3a_0
\qtq{and hence}a_1^2=4a_0\bigl(a_2+\alpha a_3\bigr).
$$
We claim that $a_1$ is a constant multiple of $a_0$. Since $a_0$ divides $a_1$
this holds once $a_0$ has only simple roots. (We need to work projectively, that is,  if $\deg a_0<2$ then there is a root at infinity of multiplicity $2-\deg a_0$.)
If $u_0$ is a double root of $a_0$ then it is also a root of $a_1$ but then
$C_S$ is singular at $x=0, u=u_0$ by (\ref{type.S.say}.3), a contradiction. 
Thus there is a $\beta\in k$ such that 
$$
a_1=2\beta a_0\qtq{and so } \beta^2a_0=\bigl(a_2+\alpha a_3\bigr).
$$
Expressing $a_2$ and substituting we get the equation
$$
a_3x^2(x-\alpha)+a_0\bigl(\beta x+1\bigr)^2.
\eqno{(\ref{2.4.specv.sec}.5)}
$$
Note further that the fiber $F_{\alpha}$ is 
isomorphic to the fiber $F_0$. 
(It is the image of $F_0$ under $m_2$.)
Thus,  $m_2(F_{\alpha})\subset S$
is a multi-section of $\pi$, except when the  polynomial  given by the substitution
$x=x'+\alpha$, namely 
$$
a_3(x'+\alpha)^2x'+a_0\bigl(\beta x'+\alpha\beta+1\bigr)^2,
\eqno{(\ref{2.4.specv.sec}.6)}
$$
 has the same form as (\ref{2.4.specv.sec}.5).
Rearranging by powers of $x'$ we get that 
$$
a'_0=a_0(\alpha\beta+1\bigr)^2
\qtq{and} a'_1=a_3\alpha^2+2a_0\beta(\alpha\beta+1\bigr).
$$
Thus $a'_0\mid a'_1$ holds iff $a_0\mid a_3$, unless $\alpha=0$. 
If $a_0\mid a_3$ then
the polynomial (\ref{2.4.specv.sec}.5) is reducible and $C_S$ is singular, so this can not happen.

Finally, if  $\alpha=0$,
we are down to just 2 exceptional types
corresponding to $\beta\neq 0$ and  $\beta=0$:
\begin{enumerate}\setcounter{enumi}{6}
\item $\bigl(z^2=a_3(u)x^3+a_0(u)(x-1)^2\bigr)$. In this case $F_1$
is another smooth fiber and  $m_2(F_1)\subset S$
is a multi-section of $\pi_1$.
\item $\bigl(z^2=a_3(u)x^3+a_0(u)\bigr)$.
 There are no other obvious smooth fibers.
\end{enumerate}
 In both  cases
$F_0$ consists of 3-torsion points and every multiplication map  $m_r$
sends $F_0$ to itself. \qed

\begin{cor} \label{unirat.deg2.cor}
Let  $S$ be  a degree 1, minimal conic bundle of Type S
over $k$. 
\begin{enumerate}
\item If $k$ is infinite then 
  there is a degree 2 field extension $K/k$ such that
 $S_K$ is unirational.
\item If  $k$ is finite and $|k|\geq 19$ then $S$ is unirational.
\end{enumerate}
\end{cor}

Proof.  Since an elliptic curve has at most 11 non-trivial  2- or 3-torsion
 points,
$\pi:S\to \p^1$ has at most 11 fibers that consist of  2- or  3-torsion points.
Pick any point $v\in \p^1(k)$ such that the fiber $F_v$ over $v$ is smooth
and does not  consist of  2- or  3-torsion points.
 Since $F_v$ is a conic,
it has a point in some degree 2 field extension $K/k$. 
Then Proposition \ref{2.4.specv.sec}
produces a rational multi-section  $\sigma:\p^1_K\to S_K$,
hence $S_K$ is unirational. 

If $k$ is finite then $F_v$ has a $k$-point, so we get unirationality over $k$.
\qed

\section{Birational maps of conic bundles}

\begin{say}[Elementary transformations]\label{el.transf.say}
 Let $\pi:S\to B$ be a
conic bundle over an algebraically closed field, $F\subset S$ a smooth fiber
and $q\in F$ a   point. The elementary transformation of $S$
centered at $q$ is obtained by first blowing up $q$ and then
contracting the birational transform of $F$. We get another
 conic bundle  of the same degree and a birational map  $\rho_q: S\map S_q$.

More generally, let $Q\subset S$ be a finite collection of
 points, each contained in a smooth fiber such that
every fiber contains at most 1 point. Let $F_Q$ be the union of all fibers that
contain a point of $Q$. We can then  blow up $Q$ 
to get $\alpha_Q:T_Q\to S$ with exceptional curve $E_Q$.
Let $F'_Q\subset T_Q$ denote the birational transform of $F_Q$.
We can next
contract  $F'_Q$ to get $\beta_Q: T_Q\to S_Q$ and the composite
$\rho_Q: S\map S_Q$.  These maps and surfaces fit into a diagram
$$
\begin{array}{c}
 T_Q\\
\alpha_Q \swarrow \qquad\searrow \beta_Q \\ 
 S \ \quad\stackrel{\rho_Q}{\map} \ \quad S_Q
\end{array}
\eqno{(\ref{el.transf.say}.1)}
$$
If $S$ and $Q$ are 
defined over $k$ then  so is
$\rho_Q: S\map S_Q$.

In order to compute the birational transform of $-K_{S_Q}$ note that
$$
K_{T_Q}\sim \alpha_Q^*K_S+E_Q\sim \beta_Q^*K_{S_Q}+F'_Q
\qtq{and}  F'_Q\sim \alpha_Q^*F_Q-E_Q.
$$
These together imply that
$$
\beta_Q^*(-K_{S_Q})\sim \alpha_Q^*(-K_S+F_Q)-2E_Q.
\eqno{(\ref{el.transf.say}.2)}
$$
Equivalently, if the linear system $|-K_{S_Q}|$
 is not empty then
$$
\rho_Q^*|-K_{S_Q}|=|-K_S+F_Q|(-2Q).
\eqno{(\ref{el.transf.say}.3)}
$$
That is, the pull-back of $|-K_{S_Q}| $ by $\rho_Q$
consists of those curves in  $|-K_S+F_Q| $
that have multiplicity $\geq 2$  at each point of $Q$. 

Consider now the very special case when there is 
  an irreducible curve $C\in |-K_S|$  and $Q$ consists of smooth points
of $C$. Then 
$$
C+F_Q\in    |-K_S+F_Q|(-2Q)
\eqno{(\ref{el.transf.say}.4)}
$$
and $C_Q:=\rho_Q(C)\in |-K_{S_Q}| $.
In particular, $S_Q$ is also weak del Pezzo by (\ref{wdP.all.say}.4). 
Note that the line bundle $\o_{S}(C+F_Q) $ does not depend on the
choice of $C$ and $Q$ since it is isomorphic to
$\o_{S}(-K_S+nF) $ where $n:=|Q|$. Thus, using  (\ref{el.transf.say}.3)
and the natural isomorphism  $C\cong C_Q$
we see that
$$
\o_{S_Q}(C_Q)|_{C_Q}\cong \o_{S}(C+F_Q)|_{C}(-2Q)\cong 
\o_{S}(-K_S+nF)|_{C}\otimes\o_C(-2Q).
\eqno{(\ref{el.transf.say}.5)}
$$
\end{say}

As a first application we describe
how to transform  conic bundles of Types S or E into a
conic bundle of Type G. 


\begin{prop} \label{typeS.bir.typeG}
Let  $S$ be a degree 1, minimal  conic bundle of Type S or E
over a field $k$ with $\chr k\neq 2$ and $|k|\geq 53$.
Then $S$ is birational to a conic bundle of Type G iff it has
at least one  $k$-point that lies on a smooth fiber.
Thus every such $S$ has a degree 8 rational  multi-section.
\end{prop}

Proof. We start with Type E. Let $C\subset S$ be the rational double
section and $F$ a smooth fiber with a $k$-point. Pick
$q\in F(k)$ that is not on $C$. We get an elementary
transformation $\rho_q:S\map S_q$.  Note that $C\in |-K_S-2F|$ by
(\ref{exceptional.cases.prop}.2), 
hence $C+2F\in |-K_S|(-2q)$ and so $C_q\in |-K_{S_q}|$. 
 Thus $S_q$ is weak del Pezzo by (\ref{wdP.all.say}.4).

Assume next that $S$ is of type S.
Pick a $k$-point $q$ lying on a  smooth fiber  $F$.
Note that $F\cong \p^1_k$ since $F$ is a conic with a $k$ point.

Let $C\in |-K_S|$ be the unique curve passing through $q$. 
It is better to think of $C$ as determined by the point pair
$F\cap C=\{q,q'\}$.

The restriction of  $|-K_S|$  to  $F$ is a degree 2 pencil on $F$, hence
it has at most 2 ramification points since $\chr k\neq 2$.
Thus  $q\neq q'$ holds for all but 2 pairs $\{q,q'\}$.
Since the elliptic pencil has at most 12 singular fibers,
$C$ is  geometrically irreducible and  $q\neq q'$
 for all but 26 points  $q\in F(k)$.

We are in the case considered in (\ref{el.transf.say}.4), thus
the map $\rho_q$ transforms $C$ into a curve $ C\cong C_q\in |-K_{S_q}|$
and, by   (\ref{el.transf.say}.5), 
$$
 \o_{S_q}(C_q)|_{C_q}\cong \o_C(3p^*-2q).
$$
(Here we used that $S$ is of Type S, thus $F|_C\sim 2p^*$.) 
Therefore the base point $p_q^*$ of  $S_q$ is the
 unique point on $ C\cong C_q$ such that
$p^*_q\sim 3p^*-2q$. 
Hence $S_q$ is of Type G iff $p_q^*$
lies on a smooth fiber.

If $p_q^*=p^*$ then $2p^*\sim 2q$ and so  $C$ is either singular
at $q$ or is tangent to
$F$ at $q$. This  was excluded. 

If $p_q^*$ lies on  a singular fiber, then it is the
unique singular point of the fiber. Since any two curves in $|-K_S|$ meet
only at $p^*$, in this case $p_q^*$ uniquely determines  
$\{q,q'\}$. This gives another 24 possible exceptions. 
Thus, with  at most 50 exceptions, $p_q^*$
lies on a smooth fiber of $\pi$ and   $S_q$ is of Type G. 
The last assertion then follows from Proposition \ref{type.G.unirtl}.\qed

\section{Unirationality for Type S; second construction}\label{sec.unir.S.2}

Let $T\subset \p^3$ be a cubic surface. Given 2 points $q,q'$ on $T$
the line connecting them usually intersects $T$ in a unique 3rd point
$p(q,q')$.
This gives  a rational map  $\sym^2(T)\map T$ that is very useful
in studying rational points on $T$; see for instance
\cite{k-looking}. Note that one can also think of
$p(q,q')$ as the base point of the linear system  $|-K_T|(-q-q')$.

Generalizing the latter point of view gives similar maps
for degree 1 conic bundles. This will complete the proof of 
Theorem \ref{main.thm}.

\begin{const}\label{deg.1.construction}
 Let $S$ be a degree 1, weak del~Pezzo conic bundle.
The linear system $|-K_S+nF| $ has dimension  $3n+1$ and
self-intersection $4n+1$. 

Therefore, given any set of $n$ points $Q=\{q_1,\dots, q_n\}$,
  the linear system  
$$|-K_S+nF|(-2Q)
$$
has dimension at least 1 and there is a---possibly empty---open 
subset $U\subset S^n$  such that 
$|-K_S+nF|(-2Q)$ has the following properties for every
$(q_1,\dots, q_n)\in U$.
\begin{enumerate}
\item  The  dimension is 1.
\item The general member is  geometrically irreducible and, 
after blowing up $Q$, it has arithmetic genus 1.
\item The base points consist of $Q$  (with multiplicity 2)
and one more point; we denote it by   $p^*_Q$.
\end{enumerate}
 We will check that $U\neq \emptyset$.
If this holds then $Q\mapsto p^*_Q$  defines a morphism  $U\to S$
which descends to a rational map
$\Phi_n:\sym^n(S)\map S$.
\end{const}

\begin{prop}\label{Phi.defined.dom.prop} Let $k$ be a field and 
  $S$  a degree 1, weak del~Pezzo conic bundle.
 Then the above  $\Phi_n:\sym^n(S)\map S$
is defined and dominant.
\end{prop}

Proof. It is sufficient to check these claims over the algebraic
closure of $k$. Thus we may assume to start with that
$k$ is algebraically closed. 
We will exhibit  a rather special set of points $Q=\{q_1,\dots, q_n\}$
such that $(q_1,\dots, q_n)\in U$.

Pick an irreducible curve $C\in |-K_S|$ and smooth points
$q_1,\dots, q_n\in C$ such that each $q_i$ is contained in a
smooth fiber $F_i$ of $\pi$ and the $F_i$ are all distinct.
Then $S_Q$ is another degree 1, weak del~Pezzo conic bundle
hence $|-K_{S_Q}|$ has dimension 1 and
its general member is  geometrically irreducible and
has arithmetic genus 1. Thus, by (\ref{el.transf.say}.3), 
$$
\rho_Q^*|-K_{S_Q}|=|-K_S+F_Q|(-2Q)
$$
also has dimension 1,
its general member is  geometrically irreducible and,
after blowing up the $q_i$,
has arithmetic genus 1. Thus the extra base point $p_Q^*$
is the preimage (under $\rho_Q$) of the
base point of $|-K_{S_Q}| $ (which we also usually denote by $p_Q^*$).
Furthermore, (\ref{el.transf.say}.5) says that
$$
\o_{C_Q}(p_Q^*)\cong\o_{S_Q}(C_Q)|_{C_Q}\cong 
\o_{S}(-K_S+nF)|_{C}\otimes\o_C(-2Q).
$$
Thus $\Phi_n(q_1,\dots, q_n)$ is the unique point on $C$ such that
$$
\Phi_n(q_1,\dots, q_n)\sim p^*+n(F|_C)-2\tsum_iq_i.
$$
By varying the points $q_i$ inside $C$ we see that $\Phi(U)\cap C$ 
is dense in   $C$ and by varying $C$ we conclude  that 
 $\Phi(U)$ is dense in $S$. \qed

\begin{rem} Note that while taking the points $q_i$
on the same member of $|-K_S|$ makes the geometric computation easy,
it does not give   new $k$-points unless $C$ already has
nontrivial $k$-points. 
\end{rem}

\begin{prop} \label{symm.power.unirat.prop}
Let $X$ be a geometrically irreducible $k$-variety and $K/k$
 a  field extension  of degree $n$.   
If $X$ is unirational over $K$ then
$\sym^n(X)$ is unirational over $k$.
\end{prop}

Proof. Since $X_{K}$ is unirational, there is a
dominant rational map  $\phi_K:\a^r\map X$ defined over $K$. 
Taking the Weil restriction (see, for instance  \cite[Sec.7.6]{blr})
gives a dominant rational map
$$
\Re_{K/k}\bigl(\phi_K\bigr):\a^{rn}_k\cong \Re_{K/k}\bigl(\a^r\bigr)\map \Re_{K/k}(X),
$$
 defined over $k$. Composing with the natural map
$\Re_{K/k}(X)\to \sym^n(X)$ shows that $\sym^n(X)$ is unirational over $k$.
\qed

\begin{cor}\label{deg1.wDP.unir.bifbd.cor}
 Let $k$ be a field and
$S$  a degree 1, weak del~Pezzo conic bundle over $k$.
Then $S$ is unirational over $k$.
\end{cor}

Proof. Since $S$ is geometrically rational, it is
unirational over a finite field extension $K/k$. 
By Proposition \ref{symm.power.unirat.prop}
 we conclude that $\sym^n(X)$ is unirational over $k$
where $n=\deg (K/k)$. Thus $X$ is unirational by 
Proposition \ref{Phi.defined.dom.prop}.\qed
\medskip

\begin{rem}
Note that the above argument proves the unirationality
part of Theorem \ref{main.thm}
without using any of the earlier  constructions.
However, a direct application of this argument
gives a rather large value for $n$. The Galois group
of $\bar k/k$ acts on the 14 $(-1)$-curves
in the  7 singular fibers. Thus all of these  $(-1)$-curves
are defined over a Galois extension $K/k$ whose 
degree divides  $7!2^7$. 
Over $K$ we can contract 7 of these curves to see that
$S_K$ is birational to a $\p^1$-bundle over $\p^1$, hence rational.
Thus 
Corollary \ref{deg1.wDP.unir.bifbd.cor} gives 
a dominant rational map
$$
\Phi_n:\p^{2n}\map S\qtq{where} n=7!2^7=645120.
$$
\end{rem}

\begin{rem}\label{18.8.rem} In order to get something that is more computable,
one should use   Corollary \ref{unirat.deg2.cor} 
and apply Proposition \ref{symm.power.unirat.prop}
with $n=2$. We have not been able to write an explicit formula
for $\Phi_2$ but
the proof gives the
following description of the 8-fold section
for the equation $(**)$.

First some notation. Let  $t_0, t_0'$ be the roots of $a_3(t)=0$
and set $p:=\bigl( (1{:}0),(t_0{:}1)\bigr)$
and $p':=\bigl( (1{:}0),(t'_0{:}1)\bigr)$.
Given polynomials  $f(u),g(u), F(u), G(u)\in K[u]$ let 
 $C(f,g,F,G)\subset\p^1_x\times \p^1_t$  denote the 
(closure of the) image of the map
$$
u\mapsto  \bigl(f(u)/g(u), F(u)/G(u)\bigr).
$$

{\it Step 1.} Pick $x,y\in K$ at random. Solving $(**)$ for $t$
we get a conjugate pair of points  $r,r'\in S$.

{\it Step 2.} Using the elliptic curves $(*)$, compute
$r_2\sim -2r$ and  $r'_2\sim -2r'$. Let $q,q'\in \p^1_x\times \p^1_t$
denote the projections of $r_2,r'_2$.

{\it Step 3.} We are now looking for the unique member $C^*$ 
of the rational pencil
$|-4K_S+7F|(-8q-8q')$ that passes through $p^*$. 
Choosing a birational map $\p^1_u\to C^*$ such that $\infty\mapsto p^*$,
in the representation
of Paragraph \ref{type.S.say} the map  $\p^1_u\to C^*\into S$
 is given as follows. 
\begin{enumerate}
\item[(i)] There are
polynomials  $f,g, F, G, H$ of degrees 8,6,17,17,26
 such that
\item[(ii)] $C(f,g,F,G)$ has multiplicity 8 at $q,q'$,
\item[(iii)] $C(f,g,F,G)$ has multiplicity 3 at $p,p'$ and
\item[(iv)] $a_3(F/G)(f/g)^3+\cdots+ a_0(F/G)=\bigl(H/(gG)\bigr)^2$.
\end{enumerate}
(Explanation: We know that $C^*$ is a degree 8 multi-section, thus
$f(u)/g(u)$ is a rational function of degree 8. We chose 
 $\p^1_u\to C^*$ such that $\infty\mapsto p^*$.
Since $p^*$ is a singular point of the fiber at infinity, this
corresponds to a double pole at infinity, thus $\deg g=6$.
Furthermore, $C^*\in |-4K_S+7F|$ and so $(-K_S\cdot C^*)=18$.
Blowing up $p^*$ reduces the intersection number by 1, hence
$\deg F, \deg G\leq 17$.  The expression on the left side of (iv) 
has denominator
$g^3G^2$ but this can be a square only if a $g$ factor gets canceled.
This gives $\deg H=26$.)

{\it Step 4.}  The 8-fold section is the image of the map
$$
u\mapsto\left(\frac{f(u)}{g(u)}, \frac{H(u)}{g(u)G(u)},
\frac{F(u)}{G(u)}\right).
$$
We can thus obtain many rational points using the above map.

{\it Step 5.} Once we have a rational point $q=(x_0, y_0, t_0)\in S$, 
we also have a simpler  degree 8 rational multi-section.
In this case we
 need to find polynomials  $f,g, F, G, H$ of degrees 8,6,9,9,18
 such that $C(f,g,F,G)$ has multiplicity 8 at $q$ and
3 at $p,p'$. The rest is the same as above.

{\it Step 6.} A subtlety with this method is that these formulas work only for
points in general position and it is difficult to describe  these conditions
explicitly.  Another problem is that $\aut(\p^1,\infty)$ acts on all solutions
and it is not clear how to rigidify the set-up.

\end{rem}

\section{Unirationality of higher degree conic bundles}
\label{deg2.section}

We complete the proof of Corollary \ref{main.thm.2} by proving that
 minimal  conic bundles $\pi:S\to \p^1$ of degree $\geq 2$
are unirational over $k$ provided they have a $k$-point.

As we noted after  Corollary \ref{main.thm.2}, the only possibly new case
is the unirationality of degree 2, minimal  conic bundles.
We  reduce it to the already established degree 1 case. 
One could use the methods of the previous sections
to get a more direct argument but some case analysis
 would still remain.

As before, we consider only fields whose characteristic is
different from 2. The main reason is that many proofs involve
conics or involutions and in this way we avoid inseparability problems.
There would be  further difficulties with very small fields, especially
with $\f_2$.

\begin{prop} \label{deg2.unirat.prop}
Let $k$ be a field and $\pi:S\to \p^1$  a degree 2, minimal  conic bundle 
over $k$ with a $k$-point $p\in S(k)$.  Then $S$ is unirational over $k$.
\end{prop}

Proof. By Lemma \ref{deg2.birat.wDP} it is enough to
consider the case when   $S$ is  weak del~Pezzo. 

If $p$ lies on a smooth fiber of $\pi$ then the blow-up
$B_pS$ is a  degree 1  conic bundle 
over $k$ hence unirational by  Corollary \ref{deg1.wDP.unir.bifbd.cor}.

If $p$ lies on a singular fiber $F_p$, then $B_pS$ is not a conic bundle.
We can however contract the birational transform of $F_p$
to get a degree 1  conic bundle 
over $k$, albeit with a conjugate pair of $A_1$-singularities.
It is easy to see that  Corollary \ref{deg1.wDP.unir.bifbd.cor}
works for degree 1  conic bundles with Du~Val singularities.

 It is, however, cleaner to follow Construction \ref{deg.1.construction}
 directly. 
Given  $Q=\{q_1,\dots, q_n\}$ we consider the linear system 
$|-K_S+nF|(-2Q-p)$.
As in  Proposition \ref{Phi.defined.dom.prop}   we show that,
for general $Q$, $|-K_S+nF|(-2Q-p)$ has a unique
 unassigned base point $p_Q^*$ and it defines a dominant rational map
$\Phi_n:\sym^n(S)\map S$. Thus $S$ is unirational 
by Proposition \ref{symm.power.unirat.prop}. \qed

\begin{lem} \label{deg2.birat.wDP}
Let $k$ be a field and $\pi:S\to \p^1$  a degree 2, minimal  conic bundle 
over $k$.  Then $S$ is birational 
to a  degree 2, minimal, weak del~Pezzo conic bundle.
\end{lem}

Proof. There is nothing to do if $S$ is  weak del~Pezzo.
Otherwise,  we are in the exceptional case (\ref{exceptional.cases.prop}.2), 
hence
$S$ contains a conjugate pair of disjoint $(-3)$-curves  $A,A'$ 
and $|-K_S|=A+A'+|2F|$.

Pick a fiber $F_1$ defined over $k$. Then $A\cap F_1$ determines
a degree 2 field extension $k'/k$. Assume  now that there is a smooth fiber
$F_2$ that is defined over $k'$ but not over $k$. Then
 $A\cap F_2$ is a $k'$-point, thus $F_2$ is rational over $k'$.
Thus it has a $k'$-point $r$  not on $A+A'$; let $\bar r$
denote its conjugate. Since $F_2$ is  not defined  over $k$,
$r,\bar r$ lie on different smooth fibers.

 The elementary transformation centered at $r+\bar r$
gives  $S\map S_r$. The birational transform of $A+A'$
is $A_q+A'_q\in |-K_{S_r}|$, a conjugate pair of  $(-1)$-curves
meeting at 2 points. Thus $S_r$ is a weak del~Pezzo surface by
(\ref{wdP.all.say}.4).

One can always find the required $F_2$, save when $k=\f_3$
and the 6 singular fibers lie over the 6  points of 
$\f_9\setminus \f_3$. In this case the fibers over the $\f_3$ points
are smooth and rational. We can thus pick 2 $\f_3$-points
in 2 different smooth fibers to play the role of $r,\bar r$. 

(Actually, there is no minimal conic bundle where  the 6 singular fibers lie over the 6  points of 
$\f_9\setminus \f_3$. Indeed, these points form 3 conjugacy classes,
but class field theory (see for instance \cite[Sec.XIII.6]{MR0234930})
implies that,  over a finite field, the  singular fibers form an even number of conjugacy classes.)

\qed


\def\cprime{$'$} \def\cprime{$'$} \def\cprime{$'$} \def\cprime{$'$}
  \def\cprime{$'$} \def\cprime{$'$} \def\dbar{\leavevmode\hbox to
  0pt{\hskip.2ex \accent"16\hss}d} \def\cprime{$'$} \def\cprime{$'$}
  \def\polhk#1{\setbox0=\hbox{#1}{\ooalign{\hidewidth
  \lower1.5ex\hbox{`}\hidewidth\crcr\unhbox0}}} \def\cprime{$'$}
  \def\cprime{$'$} \def\cprime{$'$} \def\cprime{$'$}
  \def\polhk#1{\setbox0=\hbox{#1}{\ooalign{\hidewidth
  \lower1.5ex\hbox{`}\hidewidth\crcr\unhbox0}}} \def\cdprime{$''$}
  \def\cprime{$'$} \def\cprime{$'$} \def\cprime{$'$} \def\cprime{$'$}
\providecommand{\bysame}{\leavevmode\hbox to3em{\hrulefill}\thinspace}
\providecommand{\MR}{\relax\ifhmode\unskip\space\fi MR }
\providecommand{\MRhref}[2]{%
  \href{http://www.ams.org/mathscinet-getitem?mr=#1}{#2}
}
\providecommand{\href}[2]{#2}

\vskip1cm

\noindent Princeton University, Princeton NJ 08544-1000

{\begin{verbatim}kollar@math.princeton.edu\end{verbatim}}

\medskip

\noindent Dipartimento di Matematica e Informatica, Universit\`a di Ferrara, 

\noindent Via
Machiavelli 35, 44100 Ferrara Italia

{\begin{verbatim}mll@unife.it\end{verbatim}}

\end{document}